\documentclass[12pt]{article} 

\usepackage[left= 0.8in, textwidth=6.9in, textheight=9.0in, top=1.0in]{geometry}

\usepackage{graphicx}
\usepackage{amsmath,amssymb,verbatim}
\usepackage{color}  
\usepackage{bm}                                  

\newcommand{\MT}{\left[ \begin{array}{rrrrrrrrrrrrrrrrrrrr}}
\newcommand{\EM}{\end{array}\right]}

\def\Fr{\ds \frac}
\def\ds{\displaystyle}
\def\dsum{\ds\sum}
\def\subEQ{\begin{subequations}}
\def\subEE{\end{subequations}}

\title{\LARGE \bf
Algorithms of Data Development For Deep Learning and Feedback Design
\thanks{This work was supported in part by U.S. Naval Research Laboratory - Monterey, CA}
}

\author{Wei Kang\thanks{Department of Applied Mathematics, Naval Postgraduate School, Monterey, CA, USA; wkang@nps.edu}
\and Qi Gong\thanks{Department of Applied Mathematics,University of California at Santa Cruz, Santa Cruz, CA, USA}
\and Tenavi Nakamura-Zimmerer\thanks{Department of Applied Mathematics,University of California at Santa Cruz, Santa Cruz, CA, USA}}

\begin{document}

\maketitle

\begin{abstract}
Recent research reveals that deep learning is an effective way of solving high dimensional Hamilton-Jacobi-Bellman equations. The resulting feedback control law in the form of a neural network is computationally efficient for real-time applications of optimal control. A critical part of this design method is to generate data for training the neural network and validating its accuracy. In this paper, we provide a survey of existing algorithms that can be used to generate data. All the algorithms surveyed in this paper are causality-free, i.e., the solution at a point is computed without using the value of the function at any other points. At the end of the paper, an illustrative example of optimal feedback design using deep learning is given.
\end{abstract}

\section{Introduction}
In the design of feedback controllers for dynamical systems, a critical part is to find a mathematical feedback law that is guaranteed to meet the performance requirement for the system model. The performance of feedback controls includes, but is not limited to, minimizing a cost function, stabilization, tracking, synchronization, etc. After decades of active research on control theory with tremendous progress, there exists a huge literature of feedback design methodologies. For linear systems, there are a well developed theory and commercially available computational tools, such as MATLAB toolboxes, that implement the linear theory for practical applications. For nonlinear control systems, rigorous theory and design methods have been developed. However, their applications lag behind in many areas. Lacking effective computational algorithms for nonlinear feedback design is a main bottleneck. For instance, ideally one would like to solve the Hamilton-Jacobi-Bellman (HJB) equation so that a simple control law can be derived for optimal control. However, solving the HJB equation is a problem that suffers the curse-of-dimensionality. For systems with even a moderate dimension such as $n\geq 4$, finding a numerical solution for the HJB equation is extremely difficult, if not impossible. The curse-of-dimensionality affects many areas of nonlinear feedback design, such as differential games, reachable sets, the PDE (or FBI equation) for output regulation, stochastic control, etc. 

Some recent publications reveal that  neural networks can be used as an effective tool to overcome the curse-of-dimensionality. In \cite{han,naka,izzo}, neural network approximations of optimal controls are found for examples with high dimensions from $n=6$ to $n=100$. In these examples, a large number of numerical solutions, or data, are first generated for the problem of optimal control. Then a neural network is trained through supervised learning based on the data. These examples show that a neural network is not sensitive to the increase of the dimension, but requires a lot of data. The methods in \cite{han,naka,izzo} share something in common. They are model-based and data-driven, i.e., the control law is designed through the training of a neural network based on data that is generated using a numerical model. In this paper, we give a survey of some existing algorithms of open-loop optimal control that can be used to generate data for the purpose of training a neural network. The survey does not include the results, such as \cite{huang} for output regulation, \cite{raissi} for general PDEs and \cite{bittracher} for multiscale stochastic systems, that are not focused on optimal control although they share similar ideas of deep learning for dynamic systems. 

For background information, we briefly introduce the key steps in the process of training a neural network to approximate an optimal control. Following \cite{naka}, consider the following problem
\begin{equation}
\label{eq: OCP}
\left \{
\begin{array}{cl}
\underset{\bm u \in \mathcal U}{\text{minimize}} &  \displaystyle \int_{t_0}^{t_f} L (t, \bm x, \bm u) dt + \psi (\bm x (t_f)), \\
\text{subject to} & \dot {\bm x} (t) = \bm f (t, \bm x, \bm u) , \\
	& \bm x (t_0) = \bm x_0 .
\end{array}
\right .
\end{equation}
Here $\bm x (t) : [t_0, t_f] \to \mathcal X \subseteq \mathbb R^n$ is the state, $\bm u (t, \bm x) : [0, t_f] \times \mathcal X \to \mathcal U \subseteq  \mathbb R^m$ is the control, $\bm f (t, \bm x, \bm u) : [0, t_f] \times \mathcal X \times \mathcal U \to \mathbb R^n$ is a Lipschitz continuous vector field, $\psi (\bm x (t_f)) : \mathcal X \to \mathbb R$ is the terminal cost, and $L (t, \bm x, \bm u) : [0, t_f] \times \mathcal X \times \mathcal U \to \mathbb R$ is the running cost, or the Lagrangian. The optimal cost, as a function of $(t_0,\bm x_0)$, is called the value function, which is denoted by $V(t_0,\bm x_0)$ or simply $V(t,\bm x)$.  The following approach is from \cite{naka}. An illustrative example is given in Section \ref{sec_example}

\begin{enumerate}
\item \textit{Initial data generation:} For supervised learning, a data set must be generated. It contains the value of $V(t,\bm x)$ at random points in a given region. A key feature desired for the algorithm is that the computation should be {\it causality-free}, i.e. the solution of $V(t_0,\bm x_0)$ is computed without using an approximated value of $V(t,\bm x)$ at nearby points. For instance, finite difference methods for solving PDEs are not causality-free (in space) because the solution is propagated over a set of grid points. The Causality-free property is important for several reasons: (1) the algorithm does not require a grid so that the computation can be applied to high dimensional problems; (2) data can be generated at targeted region for adaptive data generation; (3) the accuracy of the trained neural network can be checked empirically in a selected region; (4) data can be generated in parallel in a straightforward manner.  

\item \textit{Training:} Given this data set, a neural network is trained to approximate the value function $V(t,\bm x)$. The accuracy of the neural network can be empirically checked using a new data set. If necessary, an adaptive deep learning loop can be applied. In each round, one can check the approximation error and then expand the data set in regions where the value function is likely to be steep or complicated, and thus difficult to learn. 

\item \textit{Validation:} The training process stops when it satisfies some convergence criteria. Then, the accuracy of the trained neural network is checked on a new set of validation data computed at Monte Carlo sample points. Once again, a causality-free algorithm is needed here. 

\item \textit{Feedback control:} For real-time feedback control, causality-free algorithms may not be fast or reliable enough for most applications. However, one can compute the optimal feedback control online by evaluating the gradient of the trained neural network and applying Pontryagin's maximum principle. Notably, evaluation of the gradient is computationally cheap even for large $n$, enabling implementation in high-dimensional systems.
\end{enumerate}

Physics laws and first principle models are fundamental and critical in control system designs. Guaranteed system properties by physics laws and mathematics analysis are invaluable. These properties should be carried through the design, rather than reinventing the wheel by machine learning. In a model-based data-driven approached outlined in Steps 1 - 4, we take the advantage of existing models and design methodologies that have been developed for decades, many with guaranteed performance. The deep neural network is used focusing on the curse-of-dimensionality only, an obstacle that classical analysis or existing numerical methods have failed to overcome. Control systems designed in this way should have the performance as proved in classical and modern control theory, however be curse-of-dimensionality free.  Demonstrated in \cite{naka,han,izzo}, some advantages of the model-based data-driven approach include: the optimal feedback can be learned from data over given semi-global domains, rather than a local neighborhood of an equilibrium point; the level of accuracy of the optimal control and value function can be empirically validated; generating data using causality-free algorithms has perfect parallelism; the inherent capacity of neural networks for dealing with high-dimensional problems makes it possible to solve HJB equations that have high dimensions. 

Generating data is critical in three of the four steps shown above. A reliable, accurate and causality-free algorithm to compute $V(t,\bm x)$, and $V_{\bm x}(t,\bm x)$ in some cases, is required. Some computational algorithms for open-loop optimal control are suitable for this task. The goal of the following sections is to provide a survey of some representative algorithms that have been, or have the potential to be, used for data generation.

\section{Characteristic methods}
\label{sec_char}
Consider the problem of optimal control defined in (\ref{eq: OCP}). Let's define the Hamiltonian
\begin{equation}
\label{eq: Hamiltonian}
H (t, \bm x, \bm \lambda, \bm u) = L (t, \bm x, \bm u) + \bm \lambda^T \bm f (t, \bm x, \bm u) ,
\end{equation}
where $\bm x \in \mathbb R^n$ is the state of the control system,
$$
\dot{\bm x}=\bm f (t, \bm x, \bm u),
$$
in which $\bm u\in \mathcal U \subseteq \mathbb R^m$ is the control variable. In (\ref{eq: Hamiltonian}), $\bm \lambda \in \mathbb R^n$ is the costate and $L (t, \bm x, \bm u) $ is the Lagrangian of optimal control. The HJB equation is 
\begin{equation}
\label{HJB}
\left\{
\begin{array}{lll}
V_t(t,\bm x)+\ds\min_{\bm u \in  \mathcal U}\{ L(t,\bm x,\bm u)+V_{\bm x}^T(t,\bm x) \bm f(t,\bm x,\bm u)\}=0,\\
V (t_f, \bm x) = \psi (\bm x), 
\end{array}
\right.
\end{equation}
where $\psi(\bm x)$ represents the endpoint cost. The optimal feedback control law is
\begin{equation}
\label{eq: optimal feedback control using dVdx}
\bm u^* (t, \bm x) = \underset{\bm u \in \mathcal U}{\text{arg min}} \, H \left( t, \bm x, V_{\bm x}, \bm u \right) .
\end{equation}
The characteristics of the HJB equation follows Pontryagin's maximum principle (PMP)
\begin{equation}
\label{tpbvp}
\left\{
\begin{array}{lll}
\dot{\bm x} (t) = \Fr{\partial H}{\partial \bm \lambda} = \bm f (t, \bm x, \bm u^* (t, \bm x, \bm \lambda)) , &\bm  x(0)=\bm x_0,\\
\dot{\bm \lambda} (t) = - \Fr{\partial H}{\partial \bm x} (t, \bm x, \bm \lambda, \bm u^* (t, \bm x, \bm \lambda)) , &\bm \lambda (t_f) =\Fr{\partial \psi}{\partial \bm x}(t_f),\\
\dot v (t) =  L (t, \bm x, \bm u^* (t, \bm x, \bm \lambda)) , & v(t_f)=\psi(\bm x(t_f)).
\end{array}
\right.
\end{equation}
It is a two-point boundary value problem (TPBVP). Computational algorithms of solving a TPBVP have been studied for decades with an extensive literature. As a starting point, the four-point Lobatto IIIa formula is an algorithm that can be implemented with a controlled true error \cite{kierz}. Packages based on this algorithm exist, such as in MATLAB and in Python. As an example, the packages are applied to a TPBVP for the optimal control of a rigid-body equipped with momentum wheels in \cite{kangwilcox1,kangwilcox2,naka}. However, they have to be used with caution as explained in the following. 

\subsection{Time-marching}
\label{subsec_timemarching}
TPBVPs have been perceived as very difficult because their numerical algorithms tend to diverge. The single most important factor that affects the convergence is the initial guess. In a worst scenario,  the TPBVP solver in \cite{naka} converges at only $1\%$ of the sample pints.  However, applying a time-marching method, the convergence is improved to $98\%$. Applying more sophisticated tuning of marching steps in \cite{kangwilcox1,kangwilcox2}, $100\%$ convergence was achieved at more than $40,000$ grid points. 

In the time-marching trick, a sequence of solutions is computed that grows from an initially short time interval.
More specifically, we choose a time sequence,
$$t_0 < t_1 < t_2< \cdots< t_K = t_f,$$
in which $t_1$ is small. For the short time interval $[t_0, t_1]$, the TPBVP solver always converges using an initial guess close to the initial state. Then the resulting trajectory, $(\bm x^1(t), \bm \lambda^1(t))$  is extended over the longer time interval $[t_0, t_2]$. A simple way to extend the trajectory is with a piecewise function
$$\bm x^2_0(t)=\left\{\begin{array}{lll} \bm x^1(t), & \mbox{ if } t_0\leq t\leq t_1,\\ \bm x^1(t_1), & \mbox{ if } t_1 < t\leq t_2,  \end{array}  \right.$$
and $\bm \lambda^2_0(t)$ is similarly defined. Or one can try a linear extension 
$$\bm x^2_0(t)=\bm x^1\left(t_0+\frac{t_1-t_0}{t_2-t_0}(t-t_0)\right), \;\;\; \mbox{ for } t_0\leq t\leq t_2.$$
The trajectory over the extended interval is used as an initial guess to find $(\bm x^2(t), \bm \lambda^2(t))$, a solution of the TPBVP over $[0, t_2]$. Repeating this process until $t_K = t_f$ at which we obtain the full solution. One needs to tune the time sequence $\{t_k\}_{k=1}^K$ to achieve convergence while maintaining acceptable efficiency. The time-marching approach does not require a good initial guess. When if converges, the TPBVP solver achieves highly accurate solutions. The algorithm in \cite{kierz} even provides an estimated error. However, this method is usually slower than the neural network warm start, which is illustrated next. 

\subsection{Neural network warm start}
In \cite{naka}, a neural network warm start is used to speed up the computation and to improve the convergence. Before a neural network can be trained, we need an initial set of data. This can be generated using any causality-free algorithm, such as the time-marching method in Section \ref{subsec_timemarching}. Then we train a neural network based on the initial data set. The loss function used in \cite{naka} takes into consideration both the value functions, $V(t,\bm x)$, and the costate, $\bm \lambda (t)$. Suppose
$$V^{(i)}=V(t^{(i)}, \bm x^{(i)})$$ 
be the optimal cost, i.e. the value function, at sampling points $(t^{(i)}, \bm x^{(i)})$ for $i=1,2,\cdots, N_{d}$. If the value function is evaluated by solving the TPBVP (\ref{tpbvp}), as a byproduct the costate is also known, 
 $$\bm \lambda^{(i)}=\bm \lambda(t^{(i)}, \bm x^{(i)}),$$
where $\bm \lambda(t,\bm x)$ represents the value of costate in the solution of TPBVP with initial time $t$ and initial state $\bm x$.  Then we use a neural network to approximate $V(t,\bm x)$ by minimizing the following loss function,
$${\cal L}(\bm \theta) = \Fr{1}{N_{d}}\left(\dsum_{i=1}^{N_{d}}\left[V^{(i)} - V^{NN}(t^{(i)}, \bm x^{(i)}; \bm \theta)\right]^2+ \mu \dsum_{i=1}^{N_{d}}\left\Vert \bm \lambda^{(i)} - V_{\bm x}^{NN}(t^{(i)}, \bm x^{(i)}; \bm \theta) \right \Vert^2\right),
$$
where $\bm \theta$ is the parameter of the neural network, $\mu$ is a parameter weighing between the losses of the value function and the costate. The  trained neural network, $V^{NN}(t,\bm x)$,  provides the costate at any given point $(t,\bm x)$, which is $V^{NN}_{\bm x}(t,\bm x)$. If the size of an initial data set is small, then the neural network approximation is not necessarily accurate. However, it is good enough for the purpose of generating initial guess at any given $(t,\bm x)$ to warm start the TPBVP solver so that additional data can be generated at a much faster rate. In \cite{naka}, neural network warm start exceeds $99\%$ convergence for the rigid body optimal control problem. The initial data used to train the neural network is very small, $N_{d}=64$.

\subsection{Backward propagation}
As its name indicates, the basic idea of backward propagation is to integrate the ODEs in (\ref{tpbvp})
backward in time. In the first step, a solution of TPBVP is solved for a nominal trajectory using a nominal initial state $\bm x_0$. The second step is perturbing the final state and costate around the nominal trajectory, subject to the terminal condition. Then the ODEs in (\ref{tpbvp}) are solved backward in time to propagate a trajectory using the perturbed final state and costate value. In this way, a data set consisting of trajectories around the nominal trajectory is generated. Then, a neural network is trained by minimizing a loss function. In this approach, one avoids solving TPBVP repeatedly. Instead, the data set is generated by integrating differential equations, a task much easier than solving a TPBVP. However, the location of the sample states cannot be fully controlled. Along unstable trajectories (backward in time), integrating the ODE over a relatively long time interval can be numerically challenging. 

Backward propagation is used in \cite{izzo} for the optimal control of spacecraft making interplanetary transfers. The optimal control policy is computed for a spacecraft equipped with nuclear electric propulsion system. The goal is to transfer from the Earth to Venus orbit within about $1.376$ years. After generating $45\times 10^6$ data samples, neural networks  consisting multiple layers are used to approximate the control policy as well as the value function. Then the accuracy is validated, once again, using data generated by backward propagation.  

\section{Minimization-based methods}
Optimal control is essentially a problem of minimization (or maximization) subject to the constraint of a control system. There exist various ways of transforming them to unconstrained optimization, which can be numerically solved. 
 
 \subsection{The Hopf formula}
Consider a HJ PDE,
\begin{equation}
\label{hj}
\left\{
\begin{array}{lll}
V_t(t,\bm x) +H(V_{\bm x}(t,\bm x))=0, &\mbox{ in } (0, \infty)\times \mathbb R^n,\\
V(0,\bm x)=\psi(\bm x), & \bm x\in \mathbb R^n,
\end{array}
\right.
\end{equation}
where $H: \mathbb R^n\rightarrow \mathbb R$ is continuous and bounded from below by an affine function, $\psi: \mathbb R^n \rightarrow \mathbb R$ is convex. The Hopf formula \cite{hopf} is an expression of the solution,
\begin{equation}
\label{hopf}
V(t,\bm x)=(\psi^\ast +tH)^\ast (\bm x),
\end{equation}
where the superscript `$*$' represents the Fenchel-Legendre transform. Specifically,  $f^\ast: \mathbb R^n \rightarrow {\mathbb R} \cup \{ \infty\}$ of a function (convex, proper, lower semicontinuous) $f:  \mathbb R^n \rightarrow {\mathbb R} \cup \{ \infty\}$ is defined by 
$$f^\ast (\bm z) = \sup_{\bm x\in \mathbb R^n}\{ \bm x^T \bm z -f(\bm x)\}.$$
The solution given in (\ref{hopf}) is causality-free. In \cite{darbon}, the Hopf formula is applied to the HJ equation derived from control systems in the following form,
$$\begin{array}{lll}
\dot{\bm x}(s)=f(\bm \beta (s)),\\
\bm x(t)=\bm x,
\end{array}$$
where $\bm \beta: (-\infty, T]\rightarrow \mathbb R^n$ is the control input. The cost function is
$$J(\bm x, t;\bm \beta)=\int_t^T L(\bm\beta(s))ds+\psi(\bm x(T)),$$
where $L(\bm \beta)$ and $\psi(\bm x)$ are both scalar valued functions. The goal of optimal control design is to find a feedback that minimizes $J(\bm x, t;\bm \beta)$ using admissible control inputs. The solution can be found by solving the minimization problem
$$V(t,\bm x)=-\min_v \{ \psi^\ast(\bm v)+tH(\bm v)-\bm x^T \bm v\},$$
which is equivalent to the Hopf formula (\ref{hopf}). In \cite{darbon}, it is solved by uing the split Bregman iterative approach. Within each iteration, two minimization problems are solved numerically where Newton's method is applicable under some smoothness assumptions. In addition to optimal control, the level set method is also addressed in \cite{darbon} for the viscosity solution of the eikonal equation. 

\subsection{Minimization along characteristics}
In this approach, the value function is computed by minimizing the cost along trajectories of the Hamiltonian system (\ref{tpbvp}). Rather than a TPBVP, consider the initial value problem
\begin{equation}
\label{hamiltondyn}
\left\{
\begin{array}{lll}
\dot{\bm x} (t) = \Fr{\partial H}{\partial \bm \lambda} = \bm f (t, \bm x, \bm u^* (t, \bm x, \bm \lambda)) , &\bm  x(t_0)=\bm x_0,\\
\dot{\bm \lambda} (t) = - \Fr{\partial H}{\partial \bm x} (t, \bm x, \bm \lambda, \bm u^* (t, \bm x, \bm \lambda)) , &\bm \lambda (t_0) =\bm\lambda_0,\\
\bm u^* (t, \bm x, \bm \lambda)=\ds\arg\min_{\bm u\in {\cal U}} H(t,\bm x, \bm \lambda, \bm u).
\end{array}
\right.
\end{equation}
For fixed initial state $\bm x_0$, the cost along a characteristic is a function of $\bm \lambda_0$,
\begin{equation}
\label{cost_char}
J(t_0,\bm x_0, \bm \lambda_0)=\int_{t_0}^{t_f}  L (t, \bm x, \bm u^* (t, \bm x, \bm \lambda))dt+\psi(t_f).
\end{equation}
Then the solution, $V(t_0,\bm x_0)$, of the HJB equation (\ref{HJB}) is the minimum value of  (\ref{cost_char}) along trajectories satisfying (\ref{hamiltondyn}),
\begin{equation}
\label{eq_v}
V(t_0,\bm x_0)=\min_{\bm \lambda_0} J(t_0,\bm x_0,\bm  \lambda_0).
\end{equation}
Under convexity type of assumptions, the existence and uniqueness of solutions have been studied and proved (see, for instance,  \cite{chow,yegorov}). Similar approaches are also applicable to the HJI equation of differential games. In numerical computation, algorithms of unconstrained optimization are applicable. For instance, Powell's algorithm is used in \cite{yegorov}. In \cite{chow},  coordinate descent is used with multiple initial guesses to perform the optimization. Some examples in \cite{chow} show fast convergence that may justify real-time computation. If this is the case, then a neural network training becomes unnecessary. On the other hand, algorithms with guaranteed fast convergence for real-time optimization are still an open problem in general. 

Different from the characteristic based approach, in \cite{lin} the problem of optimal control is discretized using finite difference. Then the resulting finite dimensional optimization with constraints is transformed to an unconstrained optimization based upon Lagrangian duality. The unconstrained optimization is then solved using a splitting algorithm. This approach can be classified as a direct method, a family of computational methods based on discretization of the original problem. Some algorithms of direct methods are discussed in Section \ref{sec_direct}.

\section{Stochastic process}
Model-based deep learning for optimal control was first introduced in \cite{han2016} for stochastic systems. In this approach, the optimal control law is approximated using a neural network. The learning process is based upon a data set generated from the stochastic model of the system, rather than data sets collected from experimentation. More recently, in \cite{weinan,han} solutions of the following semilinear parabolic PDEs are approximated using a  neural network,
\begin{equation}
\label{eq_parabolic}
\left\{
\begin{array}{lll}
V_t(t,\bm x)+\Fr{1}{2}{\rm Tr}(\sigma\sigma^T{\rm Hess}_{\bm x}V)(t,\bm x)+ V_{\bm x}(t,\bm x)^T\mu(t,\bm x)+H(t,\bm x,V(t,\bm x), \sigma^T(t,x)V_{\bm x}(t,x))=0,\\
V(t_f,\bm x)=\psi(\bm x).
\end{array}
\right.
\end{equation}
As a special case, HJB equations with viscosity are PDEs in the form of (\ref{eq_parabolic}). 
In this section, $\bm x\in \mathbb R^n$, $\sigma(t,\bm x)\in \mathbb R^{n\times n}$ is a matrix valued function, $\mu(t,\bm x)$ is a vector valued function, ${\rm Hess}_{\bm x} V$ represents the Hessian of $V$ with respect to $\bm x$, ${\rm Tr}(M)$ denotes the trace of a matrix $M$, $H: (-\infty, t_f]\times \mathbb R^n\times  \mathbb R^n\times  \mathbb R^n\rightarrow  \mathbb R$ is a known function. The characteristic of the PDE is a stochastic process satisfying
\begin{equation}
\label{eq_BSDE1}
\bm x(t)=\bm x_0+\int_0^t \mu (s,\bm x(s))ds+\int_0^t \sigma(s,\bm x(s))d\bm W_s,
\end{equation}
where $\bm W_t$, $0\leq t\leq  t_f$,  is an $n$-dimensional Brownian motion. The solution of (\ref{eq_parabolic}) satisfies the following backward stochastic differential equation (BSDE) \cite{pardoux1,pardoux2},
\begin{equation}
\label{eq_BSDE2}
\left\{
\begin{array}{lll}
\begin{split}
V(t,\bm x(t))&=V(0,\bm x_0)-\int_0^t H\left(s,\bm x(s), V(s,\bm x(s)), \sigma(s,\bm x(s))^TV_{\bm x}(s,\bm x(s))\right)ds \\
&+\int_0^t V_{\bm x}(s,\bm x(s))^T\sigma (s,\bm x(s))d\bm W_s,
\end{split}\\
V(t_f,\bm x(t_f))=\psi(\bm x(t_f)).
\end{array}
\right.
\end{equation}
Discretize (\ref{eq_BSDE1}) and (\ref{eq_BSDE2}), we have
\begin{equation}
\label{eq_BSDE3}
\begin{array}{lll}
\bm x(t_{n+1})\approx \bm x(t_{n})+\mu (t_n,\bm x(t_n))(t_{n+1}-t_n)+\sigma(t_n,\bm x(t_n)) (\bm W_{t_{n+1}}-\bm W_{t_n}), & n=0,1,2,\cdots,N,
\end{array}
\end{equation}
and
\begin{align}
V(t_{n+1},\bm x(t_{n+1}))&\approx V(t_n,\bm x(t_n))-H\left(t_n,\bm x(t_n), V(t_n,\bm x(t_n)), \sigma(t_n,\bm x(t_n))^TV_{\bm x}(t_n,\bm x(t_n))\right)(t_{n+1}-t_n) \nonumber\\
&+V_{\bm x}(t_n,\bm x(t_n))^T\sigma (t_n,\bm x(t_n))(\bm W_{t_{n+1}}-\bm W_{t_n}), \label{eq_BSDE4}
\end{align}
for $ n=0,1,2,\cdots, N-1$.  In (\ref{eq_BSDE4}), the value of $V(t,\bm x)$ depends on its gradient, $V_{\bm x}(t,\bm x)$, which is unknown. 
In \cite{han}, a neural network is defined to approximate the function $x\rightarrow \sigma^TV_{\bm x}(t_n,\bm x)$ at $t=t_n$. In addition to the parameters of the neural network, the initial value and initial gradient are also treated as parameters. For the data set, the first step is to generate a set of Brownian motion paths,
$$\left\{ \{\left. \bm W^i_{t_n}\}_{n=0}^{N}\right| \;\; i=1,2,\cdots, N_s \right\},$$
where $N_s$ is the total number of samples. Then, one can compute all sample paths of the stochastic process by solving (\ref{eq_BSDE3}) (forward in time),
$$\left\{ \{\left. \bm x^i(t_n)\}_{n=0}^{N}\right| \;\; i=1,2,\cdots, N_s \right\}.$$
Integrating (\ref{eq_BSDE4}) along each sample path, one can compute $\hat V(t_n, \bm x^{(i)}(t_n))$, an approximation of the value function. Note that this approximation depends on the parameters of the neural network. 
In \cite{han}, the neural network is trained using loss functions that penalize the mean square error of the terminal condition
\begin{equation}
\label{eq_finalcond}
V(t_N, \bm x(t_N))=\psi(\bm x(t_N)).
\end{equation}
For instance, one of such loss functions is
\begin{equation}
\label{eq_loss1}
l(\theta) = \mathbb E \left(\left\{ \left. |\psi(\bm x^{(i)}(t_N)) - \hat V ( t_N,\bm x^{(i)}(t_N))|^2\right| \;\; 1\leq i\leq N_s\right\} \right),
\end{equation}
where $\theta$ represents the parameters in the neural network and the unknown initial value and initial gradient of $V(t,\bm x)$. 

In another approach \cite{raissi}, the data set is generated based upon a stochastic process, which is similar to that in \cite{han} but for coupled systems. The neural network approximates the solution,  $(t,\bm x)\rightarrow V(t,\bm x)$, over the entire time interval $[0, t_f]$. Different from (\ref{eq_loss1}), the loss function in \cite{raissi} contains two parts, the residue of the equation (\ref{eq_BSDE4}) and the error of the approximated terminal condition. In \cite{weinan,han,raissi}, neural networks are trained to solve PDEs of various dimensions up to $n=100$. In these papers,  generating data is a process of integrating (\ref{eq_BSDE3}), which is straightforward and causality-free. However, different from some TPBVP solvers for deterministic optimal control problems, the integration of (\ref{eq_BSDE4}) does not have an algorithm of error approximation. The accuracy analysis is based on the residue of the equations  (\ref{eq_BSDE4}) and/or (\ref{eq_finalcond}) at sample points.

\section{Direct methods}
\label{sec_direct}
In addition to ODEs, a control system may subject to constraints in the form of algebraic equations or inequalities, such as state constraints, control saturation, or state-control mixed constraints. A family of computational methods, so called direct methods, is particularly effective for finding optimal control with constraints. The basic idea is to first discretize the control system as well as the cost, then numerically compute the optimal trajectory using nonlinear programming. These methods do not use characteristics. Instead, the optimal control is found by directly optimizing a discretized cost function, thus the name direct method. Without a thorough review, interested readers are referred to \cite{betts,chry,don,elnagar,enright,hagar,fahroo,gong,gong2,kang1,kang} and references therein. 

For the purpose of generating data, direct methods are causality-free. As an example, in the the following we outline the basic ideas of pseudospectral (PS) optimal control. Consider the following problem
\begin{equation}
\label{Joriginal}
\min_{\bm u(\cdot )}{\cal J}[\bm x(\cdot), \bm u(\cdot)]   =  \int_{-1}^{1}F(\bm x(t), \bm u(t))\ dt + E(\bm x(-1), \bm x(1)),
\end{equation}
subject to 
\begin{equation}
\label{eq_syscons}
\left\{ \begin{array}{lll}
\dot {\bm x}(t)= \bm f(\bm x(t), \bm u(t)),\\
 \bm e(\bm x(-1), \bm x(1)) =  0,\\
 \bm h(\bm x(t),\bm u(t) ) \leq  0,
\end{array}\right.
\end{equation}
where $F: \mathbb R^{N_x}\times \mathbb R^{N_u}\rightarrow \mathbb R$,  $E: \mathbb R^{N_x}\times \mathbb R^{N_x} \rightarrow \mathbb R$, $\bm f: \mathbb R^{N_x}\times \mathbb R^{N_u} \rightarrow \mathbb R^{N_x}$, $\bm e: \mathbb R^{N_x}\times \mathbb R^{N_x} \rightarrow
\mathbb R^{N_e}$ and $\bm h: \mathbb R^{N_x} \times \mathbb R^{N_u} \to \mathbb R^{N_h}$ are continuously differentiable with respect to their arguments and their gradients are Lipschitz continuous. The problem is discretized at a time sequence 
$$t_0=-1 < t_1<t_2<\cdots < t_N=1.$$
There are various ways of chosing $t_k$. For example, the Legendre-Gauss-Lobatto (LGL) nodes (see, for instance, \cite{fahroo,gong}) are widely used in which $t_k$  are the roots of the derivative of the $N$th order Legendre polynomial. Let the discrete state and control variables be $\bar {\bm x}_k$ and $\bar{\bm u}_k$, then the problem defined in (\ref{Joriginal})-(\ref{eq_syscons}) is discretized to form a problem of nonlinear programming, 
\begin{equation}
\label{cost_disc}
\ds\min_{\bar{\bm u}} \bar J = \ds\sum_{k=0}^N F(\bar{\bm x}_k, \bar{\bm u}_k)w_k+E(\bar{\bm x}_0, \bar{\bm x}_N),
\end{equation}
subject to
\begin{equation}
\label{cons_disc}
\left\{
\begin{array}{lll}
\Vert\ds\sum_{i=0}^{N} \bar{\bm x}_i D_{ki}-\bm f(\bar{\bm x}_k, \bar{\bm u}_k)\Vert_{\infty}\leq \epsilon, & k=0,1,\cdots, N,\\
\Vert \bm e(\bar{\bm x}_0, \bar{\bm x}_N)\Vert_{\infty} \leq \epsilon,\\
\Vert \bm h(\bar{\bm x}_k, \bar{\bm u}_k)\Vert_{\infty} \leq \epsilon, & k=0,1,\cdots, N.
\end{array}
\right.
\end{equation}
In this discretization, the constraints are relaxed in which the value of $\epsilon$ is chosen to guarantee feasibility. In  (\ref{cost_disc})-(\ref{cons_disc}), $D_{kj}$ represent the elements in the differentiation matrix and $w_k$ represent the LGL weights \cite{gong}.  Problem (\ref{cost_disc})-(\ref{cons_disc}) can be solved using numerical algorithms of nonlinear programming. Commercial or free software packages are available for this purpose. The continuous time solution is approximated by
$$\bm x(t)\approx \ds\sum_{k=0}^N\bar{\bm x}_k\phi_k(t), \;\; \bm u(t)\approx \ds\sum_{k=0}^N\bar{\bm u}_k\psi_k(t),$$
where $\phi_k(t)$ is the Lagrange interpolating polynomial and $\psi_k(t)$ is any continuous function such that $\psi_k(t_j)=1$ if $k=j$ and $\psi_k(t_j)=0$ if $j\neq k$. Different choices of $\psi_k(t)$ are introduced in \cite{gong,gong2,kang,kang1}. The LGL PS method of optimal control are proved to have a high order convergence rate \cite{kang}. 

Direct methods based on Runge-Kutta or finite difference discretization are also widely used in applications.  Interested readers are referred to \cite{betts} and \cite{hagar}. 

\section{Characteristics for general PDEs}
The basic idea discussed in this paper is not limited to optimal control problems. In principal, any PDE that has a causality-free algorithm allows one to generate data for supervised learning and accuracy validation. Illustrated in Section \ref{sec_char}, finding solutions along characteristic  curves is a causality-free process, i.e., the solution can be found point-by-point without using a grid. For instance, a characteristic method of solving 1D conservation law was introduced in \cite{kang3}. Because of the causality-free property, the algorithm is able to provide accurate solution for systems with complicated shocks. In general, any quasilinear PDE, 
\begin{equation}
\ds\sum_{i=1}^na_i(x_1, \cdots, x_n,u)\Fr{\partial u}{\partial x_i}=c(x_1,\cdots,x_n,u),
\end{equation}
has a characteristic, $(\bar x(s), \bar u(s))$, defined by a system of ODEs
\begin{equation}
\left\{
\begin{array}{lll}
\Fr{d\bar x_i}{ds}=a_i(\bar x_1(s),\cdots,\bar x_n(s),\bar u(s)),\\
\Fr{d\bar u}{ds}=c(\bar x_1(s),\cdots,\bar x_n(s),\bar u(s)).
\end{array}
\right.
\end{equation}
Then the solution of the PDE satisfies
$$ u(\bar x_1(s),\cdots,\bar x_n(s))=\bar u(s).$$
A challenge of this method is that the characteristic curves may cross each other to form shocks. Also, the curves may not cover the entire region, thus forming rarefaction.  Nevertheless, if a unique solution can be defined based on  characteristics, the resulting algorithm is causality-free. It can be used to generate data to train a neural network as an approximate solution. The data can also be used to check the accuracy of the neural network solution. 

\section{An example of optimal attitude control}
\label{sec_example}
\begin{figure}[h!]
\centering
\includegraphics[width = 0.85 \textwidth]{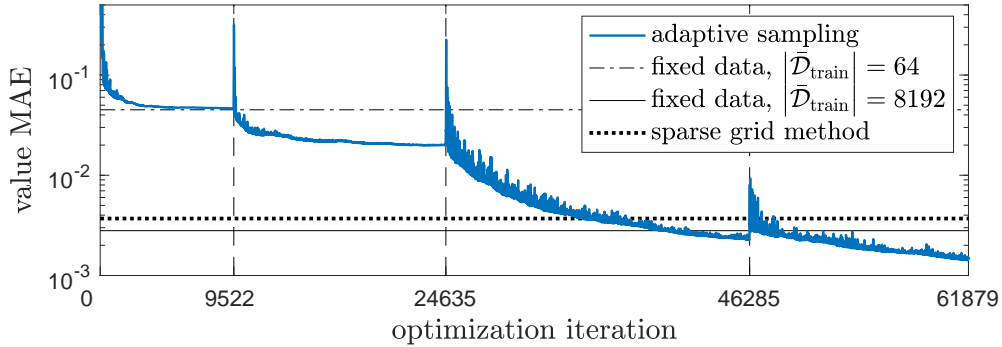}
\caption{Progress of adaptive sampling and model refinement for the rigid body problem, compared to training on fixed data sets and the sparse grid characteristics method. Vertical dashed lines show the start of new training rounds.}
\label{fig: adaptive training}
\vspace{-10pt}
\end{figure}

This is an example from \cite{naka} in which a neural network is trained using adaptive data generation for the purpose of optimal attitude control of rigid body. The actuators are three pairs of momentum wheels. 
The state variable is
$\bm x = \begin{pmatrix}
\bm v & \bm \omega
\end{pmatrix}$. Here $\bm v$ represents the Euler angles. Following the definition in  \cite{diebel},
$$
\bm v = \begin{pmatrix}
\phi & \theta & \psi
\end{pmatrix}^T ,
$$
in which $\phi$, $\theta$, and $\psi$ are the angles of rotation around a body frame $\bm e_1'$, $\bm e_2'$, and $\bm e_3'$, respectively, in the order $(1, 2, 3)$. These are also commonly called roll, pitch, and yaw. The other state variable is $\bm \omega$, which denotes the angular velocity in the body frame,
$$
\bm \omega = \begin{pmatrix}
\omega_1 & \omega_2 & \omega_3
\end{pmatrix}^T .
$$
The state dynamics are
$$
\begin{pmatrix}
\dot{\bm v} \\
\bm J \dot{\bm \omega}
\end{pmatrix}
= \begin{pmatrix}
\bm E (\bm v) \bm \omega \\
\bm S (\bm \omega) \bm R (\bm v) \bm h + \bm B \bm u
\end{pmatrix} .
$$
Here $\bm E (\bm v) , \bm S (\bm \omega), \bm R (\bm v) : \mathbb R^3 \to \mathbb R^{3 \times 3}$ are matrix-valued functions defined as
$$
\bm E (\bm v) :=
\begin{pmatrix}
1 & \sin \phi \tan \theta & \cos \phi \tan \theta \\
0 & \cos \phi & - \sin \phi \\
0 & \sin \phi / \cos \theta & \cos \phi / \cos \theta
\end{pmatrix} ,
\qquad
\bm S (\bm \omega) :=
\begin{pmatrix}
0 & \omega_3 & - \omega_2 \\
- \omega_3 & 0 & \omega_1 \\
\omega_2 & - \omega_1 & 0
\end{pmatrix} ,
$$
and
$$
\bm R (\bm v) :=
\begin{pmatrix}
\cos \theta \cos \psi & \cos \theta \sin \psi & - \sin \theta \\
\sin \phi \sin \theta \cos \psi - \cos \phi \sin \psi & \sin \phi \sin \theta \sin \psi + \cos \phi \cos \psi & \cos \theta \sin \phi \\
\cos \phi \sin \theta \cos \psi + \sin \phi \sin \psi & \cos \phi \sin \theta \sin \psi - \sin \phi \cos \psi & \cos \theta \cos \phi
\end{pmatrix} .
$$
Further, $\bm J \in \mathbb R^{3 \times 3}$ is a combination of the inertia matrices of the momentum wheels and the rigid body without wheels, $\bm h \in \mathbb R^3$ is the total constant angular momentum of the system, and $\bm B \in \mathbb R^{3 \times m}$ is a constant matrix where $m$ is the number of momentum wheels. To control the system, we apply a torque $\bm u (t, \bm v, \bm \omega) : [0, t_f] \times \mathbb R^3 \times \mathbb R^3 \to \mathbb R^m$. In this example, $m = 3$. Let
$$
\bm B = \begin{pmatrix}
1 & 1/20 & 1/10 \\
1/15 & 1 & 1/10 \\
1/10 & 1/15 & 1
\end{pmatrix} ,
\qquad
\bm J = \begin{pmatrix}
2 & 0 & 0 \\
0 & 3 & 0 \\
0 & 0 & 4
\end{pmatrix} ,
\qquad
\bm h = \begin{pmatrix}
1 \\
1 \\
1
\end{pmatrix} .
$$
The optimal control problem is
\begin{equation}
\label{eq: satellite OCP}
\left \{
\begin{array}{cl}
\underset{\bm u (\cdot)}{\text{minimize}} & \displaystyle \int_{t}^{t_f} L (\bm v, \bm \omega, \bm u) d\tau + \dfrac{W_4}{2} \Vert \bm v (t_f) \Vert^2 + \dfrac{W_5}{2} \Vert \bm \omega (t_f) \Vert^2 , \\
\text{subject to} & \dot{\bm v} = \bm E (\bm v) \bm \omega , \\
	& \bm J \dot{\bm \omega} = \bm S (\bm \omega) \bm R (\bm v) \bm h + \bm B \bm u .
\end{array}
\right .
\end{equation}
Here
$$
L (\bm v, \bm \omega, \bm u) = \frac{W_1}{2} \Vert \bm v \Vert^2 + \frac{W_2}{2} \Vert \bm \omega \Vert^2 + \frac{W_3}{2} \Vert \bm u \Vert^2,
$$
and
$$
\begin{array}{cccccc}
W_1 = 1,&
W_2 = 10,&
W_3 = \dfrac{1}{2},&
W_4 = 1,&
W_5 = 1,&
t_f = 20.
\end{array}
$$
The HJB equation associated with the optimal control has $n=6$ state variables and $m=3$ control variables. Solving the HJB equation using any numerical algorithm based on dense grids in state space is intractable because the size of the grid increases at the rate of $N^6$, where $N$ is the number of grid points in each dimension. In \cite{kangwilcox1}, a time-marching TPBVP solver is applied to compute the optimal control on a set of sparse gridpoints. In \cite{naka}, this idea is adopted to generate an initial data set from the domain
$$
\mathcal X_0 = \left \{ \left. \bm v, \bm \omega \in \mathbb R^3 \right| - \frac{\pi}{3} \leq \phi, \theta, \psi \leq \frac{\pi}{3}
\text{ and }
- \frac{\pi}{4} \leq \omega_1, \omega_2, \omega_3 \leq \frac{\pi}{4} \right \} ,
$$
This is a small data set with $N_d=64$ randomly selected initial states
$$\bm x^{(i)}=(\bm v^{(i)}, \bm \omega^{(i)}) \mbox{ for } i=1,2,\cdots, N_d.$$
Based on the data, a neural network implemented in TensorFlow \cite{abadi} is trained to approximate the value function, $V(t, \bm x)$, at $t=0$. The neural network has three hidden layers with 64 neurons in each. The optimization is achieved using the SciPy interface for the L-BFGS optimizer \cite{byrd,jones}. The loss function has two parts, 
$${\cal L}=\Fr{1}{N_d}\ds\sum_{i=1}^{N_d}\left[ V^{(i)} - V^{NN}(t^{(i)},\bm x^{(i)}; \theta)\right]^2+ \Fr{\mu}{N_d}\ds\sum_{i=1}^{N_d} \left\Vert \bm \lambda^{(i)} - V_{\bm x}^{NN}(t^{(i)},\bm x^{(i)}; \theta) \right\Vert^2,
$$
in which $\mu$ is a scalar weight. The optimization variable is $\theta$, the parameter in the neural network. The first part in the loss function penalizes the error of the neural network and the second part penalizes the error of its gradient. The loss function defined in this way takes the advantage of the fact that a TPBVP solver finds both the value of $V(t,\bm x)$ and the costate, which equals the gradient of the value function. Due to the small size of the data set, $V^{NN}$ is an inaccurate approximation of the value function. However, it is good enough to serve as the initial guess for the TPBVP solver. As a result, new data can be generated significantly faster than using time-marching. It makes adaptive data generation possible. After each training round, the location and number of additional data points are determined following a set of formulae \cite{naka}. Then a new data set is generated using a neural network warm start. In this example, a total of four training rounds are carried out in \cite{naka} with data size $N_d=64, 128, 1024$, and $ 4096$. The neural network improves its accuracy in every training round (see Figure \ref{fig: adaptive training}).

\section{Summary}
As demonstrated in the example of attitude control, causality-free algorithms generate data for not only the training of neural networks but also the validation of their accuracy. A guaranteed error upper bound is often impossible to be mathematically proved for applications of neural networks. In this case, an empirically computed approximate error provides critical information and confidence for practical applications. This is a main advantage of causality-free algorithms for the purpose of deep learning. In addition, the methods surveyed in this paper are all model-based. One has a full control of the location and amount of data to be generated, a property that is very useful for an adaptive training process. In computation, generating data using causality-free algorithms has perfect parallelism because the solution at each point is computed individually without using the function value at other points.  

\bibliographystyle{plain}
\bibliography{mybibfile}

\begin{thebibliography}{10}

\bibitem{abadi}
Mart\'{\i}n Abadi, Ashish Agarwal, Paul Barham, et~al.
\newblock {TensorFlow}: Large-scale machine learning on heterogeneous systems.
\newblock http://www.tensorflow.org/, 2015--.

\bibitem{betts}
John Betts.
\newblock {\em Practical Methods for Optimal Control Using Nonlinear
  Programming}.
\newblock SIAM, Philadelphia, 2001.

\bibitem{bittracher}
Andreas Bittracher, Stefan Klus, Boumediene Hamzi, and Christof Sch\"{u}tte.
\newblock A kernel-based method for coarse graining complex dynamical systems.
\newblock {\em arXiv:1904.08622v1}, 2019.

\bibitem{byrd}
Richard~H. Byrd, Peihang Lu, Jorge Nocedal, and Ciyou Zhu.
\newblock A limited memory algorithm for bound constrained optimization.
\newblock {\em SIAM J. Sci. Comput.}, 16:1190--1208, 1995.

\bibitem{chow}
Yat~Tin Chow, Jerome Darbon, Stanley Osher, and Wotao Yin.
\newblock Algorithm for overcoming the curse of dimensionality for
  state-dependent {Hamilton-Jacobi} equations.
\newblock {\em J. Comput. Phys.}, 387:376--409, 2019.

\bibitem{chry}
I.~Chryssoverghi, J.~Coletsos, and B.~Kokkinis.
\newblock Discretization methods for optimal control problems with state
  constraints.
\newblock {\em J. Comput. Appl. Math.}, 191:1--31, 2006.

\bibitem{darbon}
Jerome Darbon and Stanley Osher.
\newblock Algorithms for overcoming the curse of dimensionality for certain
  {Hamilton-Jacobi} equations arising in control theory and elsewhere.
\newblock {\em Res. Math. Sci.}, 3(1), 2016.

\bibitem{diebel}
James Diebel.
\newblock Representing attitude: {Euler} angles, unit quaternions, and rotation
  vectors.
\newblock
  {https://www.astro.rug.nl/software/kapteyn-beta/\_downloads/attitude.pdf},
  2006.

\bibitem{don}
A.~L. Dontchev and William~W. Hager.
\newblock The euler approximation in state constrained optimal control.
\newblock {\em Math. Comput.}, 70:173--203, 2001.

\bibitem{weinan}
Weinan E, Jiequn Han, and Arnulf Jentzen.
\newblock Deep learning-based numerical methods for high-dimensional parabolic
  partial differential equations and backward stochastic differential
  equations.
\newblock {\em Communications in Mathematics and Statistics}, 5(4):349–380,
  2017.

\bibitem{elnagar}
G.~Elnagar, M.~A. Kazemi, and M.~Razzaghi.
\newblock The pseudospectral legendre method for discretizing optimal control
  problems.
\newblock {\em IEEE Trans. Autom. Control}, 40(10):1793–1796, 1995.

\bibitem{enright}
Paul~J. Enright and Bruce~A. Conway.
\newblock Discrete approximations to optimal trajectories using direct
  transcription and nonlinear programming.
\newblock {\em J. Guid. Control Dyn.}, 15(4):994–1002, 1992.

\bibitem{fahroo}
Fariba Fahroo and I.~Michael Ross.
\newblock Costate estimation by a legendre pseudospectral method.
\newblock {\em J. Guid. Control Dyn.}, 24(2):270–277, 2001.

\bibitem{gong}
Qi~Gong, Wei Kang, and I.~Michael Ross.
\newblock A pseudospectral method for the optimal control of constrained
  feedback linearizable systems.
\newblock {\em {IEEE} Trans. Automat. Control}, 51(7):1115--1129, 2006.

\bibitem{gong2}
Qi~Gong, I.~Michael Ross, Wei Kang, and Fariba Fahroo.
\newblock Connections between the covector mapping theorem and convergence of
  pseudospectral methods for optimal control.
\newblock {\em Comput. Optim. Appl.}, 41(3):307--335, 2008.

\bibitem{hagar}
William~W. Hager.
\newblock Runge-kutta methods in optimal control and the transformed adjoint
  system.
\newblock {\em Numer. Math.}, 87(2):247–282, 2000.

\bibitem{han2016}
Jiequn Han and Weinan E.
\newblock Deep learning approximation for stochastic control problems.
\newblock {\em arXiv:1611.07422v1}, 2016.

\bibitem{han}
Jiequn Han, Arnulf Jentzen, and Weinan E.
\newblock Solving high-dimensional partial differential equations using deep
  learning.
\newblock {\em Proceedings of the National Academy of Sciences},
  115(34):8505--8510, 2018.

\bibitem{hopf}
Eberhard Hopf.
\newblock Generalized solutions of nonlinear equations of the first order.
\newblock {\em J. Math. Mech.}, 14(6):951--973, 1965.

\bibitem{izzo}
Dario Izzo, Ekin \"{O}zt\"{u}rk, and Marcus M\"{a}rtens.
\newblock Interplanetary transfers via deep representations of the optimal
  policy and/or of the value function.
\newblock {\em arXiv:1904.08809}, 2019.

\bibitem{jones}
Eric Jones, Travis Oliphant, Pearu Peterson, et~al.
\newblock {SciPy}: Open source scientific tools for {Python}.
\newblock http://www.scipy.org/, 2001--.

\bibitem{kang}
Wei Kang.
\newblock Rate of convergence for the legendre pseudospectral optimal control
  of feedback linearizable systems.
\newblock {\em J. Control Theory Appl.}, 8(4):391--405, 2010.

\bibitem{kang1}
Wei Kang, Qi~Gong, I.~Michael Ross, and Fariba Fahroo.
\newblock On the convergence of nonlinear optimal control using pseudospectral
  methods for feedback linearizable systems.
\newblock {\em International Journal of Robust and Nonlinear Control},
  17:1251--1277, 2007.

\bibitem{kangwilcox1}
Wei Kang and Lucas~C. Wilcox.
\newblock A causality free computational method for {HJB} equations with
  application to rigid body satellites.
\newblock In {\em AIAA Guidance, Navigation, and Control Conference}, AIAA
  2015-2009, Kissimmee, FL, 2015.

\bibitem{kangwilcox2}
Wei Kang and Lucas~C. Wilcox.
\newblock Mitigating the curse of dimensionality: sparse grid characteristic
  method for optimal feedback control and {HJB} equations.
\newblock {\em Comput. Optim. Appl.}, 68(2):289--315, 2017.

\bibitem{kang3}
Wei Kang and Lucas~C. Wilcox.
\newblock Solving 1{D} conservation laws using pontryagin's minimum principle.
\newblock {\em J. Sci. Comput.}, 71(1):144--165, 2017.

\bibitem{kierz}
Jacek Kierzenka and Lawrence~F. Shampine.
\newblock {BVP} solver that controls residual and error.
\newblock {\em Journal of Numerical Analysis, Industrial and Applied
  Mathematics}, 3(1-2):27--41, 2008.

\bibitem{lin}
Alex~Tong Lin, Yat~Tin Chow, and Stanley Osher.
\newblock A splitting method for overcoming the curse of dimensionality in
  {H}amilton-{J}acobi equations arising from nonlinear optimal control and
  differential games with applications to trajectory generation.
\newblock {\em arXiv:1803.01215}, 2018.

\bibitem{naka}
Tenavi Nakamura-Zimmerer, Qi~Gong, and Wei Kang.
\newblock Adaptive deep learning for high-dimensional {Hamilton-Jacobi-Bellman}
  equations.
\newblock {\em arXiv:1907.05317}, 2019.

\bibitem{pardoux1}
Etienne Pardoux and Shige Peng.
\newblock Backward stochastic differential equations and quasilinear parabolic
  partial differential equations.
\newblock In Rozovskii~B. L. and Sowers~R. B., editors, {\em Stochastic partial
  differential equations and their applications}, volume 176 of {\em Lecture
  Notes in Control and Information Sciences}, pages 200--217. Springer-Verlag
  Berlin Heidelberg, 1992.

\bibitem{pardoux2}
Etienne Pardoux and Tang Shanjian.
\newblock Forward-backward stochastic differential equations and quasilinear
  parabolic pdes.
\newblock {\em Probability Theory and Related Fields}, 114(2):123--150, 1999.

\bibitem{raissi}
Maziar Raissi, Paris Perdikaris, and George~Em Karniadakis.
\newblock Physics-informed neural networks: a deep learing framework for
  solving forward and inverse problems involving nonlinear partial differential
  equations.
\newblock {\em J. of Computational Physics}, 378:686--707, 2019.

\bibitem{huang}
Jin Wang, Jie Huang, and Stephen~S.T. Yau.
\newblock Approximate nonlinear output regulation based on the universal
  approximation theorem.
\newblock {\em International Journal of Robust and Nonlinear Control},
  10:439--456, 2000.

\bibitem{yegorov}
Ivan Yegorov and Peter~M. Dower.
\newblock Perspectives on characteristics based curse-of-dimensionality-free
  numerical approaches for solving {Hamilton-Jacobi} equations.
\newblock {\em Appl. Math. Optim.}, 2018.

\end{thebibliography}

\end{document}